\documentclass[10pt]{article}
\usepackage{fullpage}

\usepackage{amsthm,amsfonts,amssymb,mathrsfs,amsmath,latexsym,mathtools}
\usepackage{cite}
\usepackage{epsfig}
\usepackage{float}
\usepackage{graphicx}
\usepackage{hyperref}
\usepackage{comment}
\usepackage{enumerate}
\usepackage[shortlabels]{enumitem}
\usepackage{enumitem}
\usepackage{caption}
\usepackage{subcaption}
\usepackage{overpic}
\usepackage[super]{nth}
\usepackage{bbm}

\usepackage{tikz}

\theoremstyle{plain}
\newtheorem{theorem}{Theorem}[section]

\theoremstyle{definition}

\newcommand{\beq}{ \begin{equation} }
\newcommand{\eeq}{ \end{equation} }
\newcommand{\beqq}{ \begin{equation*} }
\newcommand{\eeqq}{ \end{equation*} }

\newcommand{\dd}{\mathrm{d}}
\newcommand{\ii}{\mathrm{i}}
\newcommand{\R}{\mathbb{R}}

\newcommand{\Z}{\mathbb{Z}}
\newcommand{\Prob}{\mathrm{P}}    

\newcommand{\pkmath}{\mathsf}  
\newcommand{\pkFdist}{F^{\mathrm{pKPZ}}} 
\newcommand{\pkCdisc}{C}  
\newcommand{\pkt}{\tau} 
\newcommand{\pkga}{\gamma} 
\newcommand{\pkh}{\pkmath{h}}  
\newcommand{\pkheight}{\mathbf{h}}  
\newcommand{\pK}{\pkmath{K}}  
\newcommand{\pid}{\pkmath{1}}  
\newcommand{\pspace}{\pkmath{S}}  

\newcommand{\kmath}{\mathsf}  
\newcommand{\kFdist}{F} 
\newcommand{\kt}{\tau} 
\newcommand{\kga}{\gamma} 
\newcommand{\kh}{\kmath{h}}  
\newcommand{\kheight}{\mathbf{h}} 
\newcommand{\kid}{\kmath{1}}  
\newcommand{\kK}{\kmath{K}}  
\newcommand{\ku}{u}  
\newcommand{\kv}{v}  
\newcommand{\ka}{\kmath{a}}  
\newcommand{\kb}{\kmath{b}}  
\newcommand{\kD}{\kmath{D}}  
\newcommand{\kp}{\kmath{p}}  
\newcommand{\kr}{\kmath{r}}  
\newcommand{\knt}{\kmath{t}}  
\newcommand{\kny}{\kmath{y}}  
\newcommand{\knx}{\kmath{x}}  

\newcommand{\pE}{\mathrm{E}}  
\newcommand{\pL}{\mathrm{L}}  

\newcommand{\TW}{\mathrm{TW}}  
\newcommand{\FTW}{\mathrm{F}_{\TW}}  
\newcommand{\diag}{\mathrm{diag}}  

\numberwithin{equation}{section}

\title{KPZ limit theorems}

\author{Jinho Baik\footnote{Department of Mathematics, University of Michigan,
Ann Arbor, MI, 48109, USA, \texttt{baik@umich.edu}}}

\date \today 

\begin{document}

\maketitle

\begin{abstract}
One-dimensional interacting particle systems, 1+1 random growth models, and two-dimensional directed polymers define 2d height fields.
The KPZ universality conjecture posits that an appropriately scaled height function converges to a model-independent universal random field for a large class of models. 
We survey limit theorems for a few models and discuss changes that arise in different domains. 
In particular, we present recent results on periodic domains. 
We also comment on integrable probability models, integrable differential equations, and universality. 
\end{abstract}


\section{Introduction}


The KPZ universality is concerned with, among others, one-dimensional interacting particle systems, 1+1 random growth, and two-dimensional directed polymers. 
These models define height functions $\kheight(x, t)$, two-dimensional random fields, where $x$ represents the one-dimensional spatial position and $t$ the time.  
A height function encodes the integrated current for interacting particle systems, the height for random growth models, and the free energy for directed polymer models. See Section \ref{sec:TASEP} for an example. 
The KPZ universality conjecture is that for a large class of models, the scaled height function 
\beq \label{eq:heightscaled}
	\kheight_T(\gamma, \tau)= \frac{\kheight(\gamma T^{2/3}, \tau T) - c(T)}{T^{1/3}}
\eeq
converges, up to scaling factors, to a model-independent universal 2d random field, which is called the KPZ fixed point. 
Here $c(T)$ is a non-random term determined by the macroscopic limit of the height function. 
Since the height, position, and time scale as $T^{1/3}, T^{2/3}$, and $T$, respectively, 
we say that \eqref{eq:heightscaled} is a 1:2:3 scaled height function. 
We also say that a KPZ limit theorem holds if a 1:2:3 scaled height function converges in any suitable sense for the problem at hand. 


Several physics papers \cite{Kardar-Parisi-Zhang86, Forster-Nelson-Stephen77, vanBeijeren-Kutner-Spohn85, Huse-Henley85} conjectured the 1:2:3 scale for various models in the mid-1980s. 
One of them is the paper \cite{Kardar-Parisi-Zhang86} of Kardar, Parisi, and Zhang on a nonlinear stochastic partial differential equation, now called the KPZ equation, from which the term KPZ universality is derived. 
These papers were followed by extensive research in the physics community. 
However, it remained unknown, even on a conjectural level, what the limit should be. 


The situation changed in 1999 with the publication of the paper \cite{Baik-Deift-Johansson99} by 
Baik, Deift, and Johansson, in which the authors considered the longest increasing subsequence problem of random permutations. 
This problem is equivalent to the zero-temperature free energy of a directed polymer model. 
The paper proved that the one-point distribution of an analog of the height function converges in distribution. 
See Theorem \ref{thm:Poisson} below. 
Moreover, the authors found the limiting distribution explicitly, which turned out to be the Tracy-Widom distribution from random matrix theory. 
This connection between the KPZ universality and random matrix theory was completely unexpected at that time. 
Soon after, Johansson \cite{Johansson00}  proved a similar result for another model, giving yet another example of a KPZ limit theorem. 

Exciting developments on KPZ limit theorems followed these results during the next two decades. 
For example, one-point limit theorems were extended to equal-time, multi-position distributions, multi-time distributions, and even to the 2d fields. 
Results were also generalized to several, mostly isolated, models, and algebraic underpinning of these specific models was studied. 
The 2d field limit, the KPZ fixed point, was determined, and various properties of the limit were established. 
Limit theorems were also proved for infinite space, half-infinite space, and recently finite space with periodic boundary condition. 

In this article, we give a historical overview of some KPZ limit theorems and present new results on the periodic domain case.
We start by introducing the subjects of KPZ universality, interacting particle systems, random growth, and directed polymers, in Section \ref{sec:TASEP}, focusing on one particular example. 
Then, we discuss some limit theorems on infinite spaces in Section \ref{sec:one}-\ref{sec:multitime}. 
After briefly discussing the half-infinite space case in Section \ref{sec:half}, we present new results on the periodic case in Section \ref{sec:ring}.
Section \ref{sec:dfformula} compares the formulas of the limiting multi-point distributions for infinite and periodic cases, and Section \ref{sec:integrableDE} concerns differential equations associated with distribution functions. 
We conclude the article with some comments on universality in Section \ref{sec:univ}.

The research on the KPZ universality has been developing rapidly and extensively over the last two decades.  
Hence, what is discussed in this article is only a small selection of the activities. 
The reader may benefit from other excellent survey articles such as \cite{Corwin11, Quastel12, Borodin-Gorin16} to see other aspects.

\section{TASEP, corner growth model, and exponential DLPP} \label{sec:TASEP}

This section discusses one of the most well-studied examples of one-dimensional interacting particle systems, 1+1 random growth, and two-dimensional directed polymers. 


\subsection{TASEP} 

The totally asymmetric simple exclusion process (TASEP), introduced by Spitzer \cite{Spitzer70} in 1970, is a continuous-time Markov process on $\Z$. 
At any given time, each integer site of $\Z$ is occupied by at most one particle. 
A particle moves to the adjacent site to its right after a random waiting time, but only if it is empty. The waiting time is exponentially distributed of mean $1$, and the clock starts once the neighboring site becomes vacant. 
All waiting times are independent of each other. 
Note that all moves are to the right (hence, totally asymmetric),  particles can move only one step at a time (simple), and no two particles occupy the same site at the same time (exclusion). 
Figure \ref{fig:TASEP} is an example of the configuration at a particular time. 
Black dots denote particles, and white dots mark empty sites. 
In this configuration, only three particles can move, and they do so independently of each other. 


\begin{figure}\centering
\begin{tikzpicture}[scale=0.5]
\node at (-9.2, -0.1)  {$\cdots$};
\fill (-8, 0) circle(0.2);
\fill (-7, 0) circle(0.2);
\fill (-6, 0) circle(0.2);
\draw (-5, 0) circle(0.2);
\fill (-4, 0) circle(0.2);
\draw (-3, 0) circle(0.2);
\draw (-2, 0) circle(0.2);
\fill (-1, 0) circle(0.2);
\fill (0, 0) circle(0.2);
\fill (1, 0) circle(0.2);
\draw (2, 0) circle(0.2);
\draw (3, 0) circle(0.2);
\draw (4, 0) circle(0.2);
\draw (5, 0) circle(0.2);
\node at (6.2, -0.1)  {$\cdots$};
\draw[thick, ->] (-6, 0.25) to[bend left] (-5, 0.25);
\draw[thick, ->] (-4, 0.25) to[bend left] (-3, 0.25);
\draw[thick, ->] (1, 0.25) to[bend left] (2, 0.25);
\end{tikzpicture}
\caption{TASEP}
\label{fig:TASEP}
\end{figure}

The TASEP is an example of interacting particle systems. General systems may allow, for example, left moves in addition to right moves, multi-range moves, or several particles at each site. 

One particular initial condition we focus on in this article is the step initial condition that the sites in $\Z_-\cup\{0\}$ are occupied, and all sites in $\Z_+$ are empty. 
The leftmost picture in Figure \ref{fig:corner} shows the step initial condition. 

\subsection{Corner growth model}

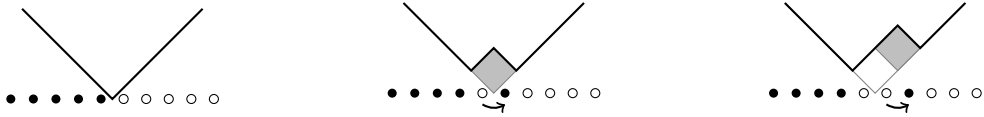
\begin{figure}
\centering
\begin{minipage}{.3\textwidth}
\centering
\begin{tikzpicture}[scale=0.3]
\fill (-4, 0) circle(0.2);
\fill (-3, 0) circle(0.2);
\fill (-2, 0) circle(0.2);
\fill (-1, 0) circle(0.2);
\fill (0, 0) circle(0.2);
\draw (1, 0) circle(0.2);
\draw (2, 0) circle(0.2);
\draw (3, 0) circle(0.2);
\draw (4, 0) circle(0.2);
\draw (5, 0) circle(0.2);
\draw[thick, -] (-3.5,4) -- (0.5,0) -- (4.5,4);
\end{tikzpicture}
\end{minipage}
\begin{minipage}{.3\textwidth}
\centering
\begin{tikzpicture}[scale=0.3]
\fill (-4, 0) circle(0.2);
\fill (-3, 0) circle(0.2);
\fill (-2, 0) circle(0.2);
\fill (-1, 0) circle(0.2);
\draw (0, 0) circle(0.2);
\fill (1, 0) circle(0.2);
\draw (2, 0) circle(0.2);
\draw (3, 0) circle(0.2);
\draw (4, 0) circle(0.2);
\draw (5, 0) circle(0.2);
\fill[lightgray] (-0.5,1) -- (0.5,2) -- (1.5,1) -- (0.5, 0) -- (-0.5, 1);
\draw[thin, gray] (-0.5,1) -- (0.5,2) -- (1.5,1) -- (0.5, 0) -- (-0.5, 1);
\draw[thick, ->] (0, -0.5) to[bend right] (1, -0.5);
\draw[thick, -] (-3.5,4)  -- (-0.5,1) -- (0.5,2) -- (1.5,1) -- (4.5,4);
\end{tikzpicture}
\end{minipage}
\begin{minipage}{.3\textwidth}
\centering
\begin{tikzpicture}[scale=0.3]
\fill (-4, 0) circle(0.2);
\fill (-3, 0) circle(0.2);
\fill (-2, 0) circle(0.2);
\fill (-1, 0) circle(0.2);
\draw (0, 0) circle(0.2);
\draw (1, 0) circle(0.2);
\fill (2, 0) circle(0.2);
\draw (3, 0) circle(0.2);
\draw (4, 0) circle(0.2);
\draw (5, 0) circle(0.2);
\draw[thin, gray] (-3.5,4)  -- (0.5,0) -- (4.5,4);
\fill[lightgray] (0.5,2) -- (1.5,3) -- (2.5, 2) -- (1.5, 1) -- (0.5, 2);
\draw[thin, gray] (0.5,2) -- (1.5,3) -- (2.5, 2) -- (1.5, 1) -- (0.5, 2);
\draw[thick, ->] (1, -0.5) to[bend right] (2, -0.5);
\draw[thick, -] (-3.5,4)  -- (-0.5,1) -- (1.5,3) -- (2.5,2) -- (4.5,4);
\end{tikzpicture}
\end{minipage}
\caption{Corner growth model}
\label{fig:corner}
\end{figure}

The configuration space for the TASEP is $\{0,1\}^{\Z}$ where $1$ represents the presence of a particle and $0$ an empty site. 
To each configuration, we can associate a zig-zag graph in $\R^2$ as in Figure \ref{fig:corner}. 
We assign a particle to a line segment of length $\sqrt{2}$ and slope $-1$, and an empty site to a line segment of the same length and slope $1$. 
Juxtaposing the line segments, we obtain a zig-zag graph as in Figure \ref{fig:corner}.
We call the graph a height function $\kheight:\R\to \R$ for the configuration of the TASEP, and it is unique up to translations. 
The leftmost picture in Figure \ref{fig:corner} is a translation of $\kheight(x)=|x|$, and it corresponds to the step initial condition. 

The TASEP induces a stochastic evolution of the height function, $\kheight(x,t)$, in which  
local valleys (corners) change to local peaks independently with rate $1$. 
The resulting 1+1 random growth process defined by the height function is called the corner growth process. 
See Figure \ref{fig:cornersim} for a simulation. 


\subsection{Exponential DLPP} \label{sec:DLPP}

Consider the two-dimensional lattice $\Z_+^2$. Let $(m,n)\in \Z_+^2$. 
An up/right (i.e., directed) path from $(1,1)$ to $(m,n)$ is a sequence $p=(p_i)_{i=1}^{m+n-1}$, where $p_i\in \Z_+^2$, $p_1=(1,1)$,  $p_{m+n-1}=(m,n)$, and $p_{i+1}-p_i\in \{(1,0), (0,1)\}$. 
The thick lines in Figure \ref{fig:ExpDLPP} are an example of a path in which we connected neighboring integer sites for visual aid. 
Let $w_s$, $s\in \Z_+^2$, be a collection of independent random variables. 
The normalized free energy of the directed polymer measure, introduced by Huse and Henley \cite{Huse-Henley85}, is 
\beqq
	F(m,n; \beta)= \frac1{\beta} \log \Big( \sum_{p} e^{\beta E(p)} \Big) \quad  \text{where $E(p)= \sum_{i=1}^{m+n-1} w_{p_i}$, }
\eeqq
the sum is over all directed paths $p$ from $(1,1)$ to $(m,n)$, and $\beta>0$ is the inverse temperature. 
The zero-temperature, $\beta=\infty$, case is called the directed last passage percolation (DLPP). 
In this case, the normalized free energy becomes 
\beqq
	L(m,n)= \max_{p} E(p), 
\eeqq
which we call the last passage time, interpreting $E(p)$ as the travel time using path $p$. 

\begin{figure}
\centering
\begin{minipage}{.4\textwidth}
\centering
\begin{tikzpicture}[scale=0.6]
\draw [->, thin][gray] (0,0) to (4.5,0);
\draw [->, thin][gray] (0,0) to (0, 3.3);
\draw [-, thin, gray] (0,0.5) to (4.3,0.5);
\draw [-, thin, gray] (0,1) to (4.3,1);
\draw [-, thin, gray] (0,1.5) to (4.3,1.5);
\draw [-, thin, gray] (0,2) to (4.3,2);
\draw [-, thin, gray] (0,2.5) to (4.3,2.5);
\draw [-, thin, gray] (0,3) to (4.3,3);
\draw [-, thin, gray] (0.5,0) to (0.5,3.3);
\draw [-, thin, gray] (1,0) to (1,3.3);
\draw [-, thin, gray] (1.5,0) to (1.5,3.3);
\draw [-, thin, gray] (2,0) to (2,3.3);
\draw [-, thin, gray] (2.5,0) to (2.5,3.3);
\draw [-, thin, gray] (3,0) to (3,3.3);
\draw [-, thin, gray] (3.5,0) to (3.5,3.3);
\draw [-, thin, gray] (4,0) to (4,3.3);
\draw [-, thick] (0,0) to (1,0) to (1,1) to (2.5, 1) to (2.5, 1.5) to (3,1.5) to (3,2) to (4,2) to (4,3);
\fill (0,0) circle(0.08);
\fill (4,3) circle(0.08);
\node at (-0.8, 0.3)  {$(1,1)$};
\node at (4.8, 3.0)  {$(m,n)$};
\end{tikzpicture}
\caption{Exponential DLPP}
\label{fig:ExpDLPP}
\end{minipage}
\qquad 
\begin{minipage}{.4\textwidth}
\centering
\begin{tikzpicture}[scale=0.6]
\fill[lightgray] (0, 3.1) to (0,2) to (1, 2) to (1, 1.5) to (1.5, 1.5) to (1.5, 1) to (2, 1) to (2, 0.5) to (3.5, 0.5) to (3.5, 0) to (4.3, 0) to (4.3, -0.2) to (-0.2, -0.2) to (-0.2, 3.1) to (0, 3.1);
\draw [->, thin][gray] (0,0) to (4.5,0);
\draw [->, thin][gray] (0,0) to (0, 3.3);
\draw [-, thin, gray] (0,0.5) to (4.3,0.5);
\draw [-, thin, gray] (0,1) to (4.3,1);
\draw [-, thin, gray] (0,1.5) to (4.3,1.5);
\draw [-, thin, gray] (0,2) to (4.3,2);
\draw [-, thin, gray] (0,2.5) to (4.3,2.5);
\draw [-, thin, gray] (0,3) to (4.3,3);
\draw [-, thin, gray] (0.5,0) to (0.5,3.3);
\draw [-, thin, gray] (1,0) to (1,3.3);
\draw [-, thin, gray] (1.5,0) to (1.5,3.3);
\draw [-, thin, gray] (2,0) to (2,3.3);
\draw [-, thin, gray] (2.5,0) to (2.5,3.3);
\draw [-, thin, gray] (3,0) to (3,3.3);
\draw [-, thin, gray] (3.5,0) to (3.5,3.3);
\draw [-, thin, gray] (4,0) to (4,3.3);
\draw [-, very thick] (0, 3.1) to (0,2) to (1, 2) to (1, 1.5) to (1.5, 1.5) to (1.5, 1) to (2, 1) to (2, 0.5) to (3.5, 0.5) to (3.5, 0) to (4.3, 0);
\end{tikzpicture}
\caption{An example of $G_t$}
\label{fig:ExpDLPPdual}
\end{minipage}
\end{figure}

For the case when $w_s\ge 0$, the DLPP is related to a random growth model. 
For $t>0$, define the subset of $\R^2$ by 
\beqq
	G_t= \bigcup_{s\in S_t} ( (0,1]^2+s)  \qquad 
	\text{where $S_t= \{(m, n)\in  \Z^2: L(m,n)\le t\}$}
\eeqq
and we set $L(m,n)=0$ if $m\le 0$ or $n\le 0$. See Figure \ref{fig:ExpDLPPdual}. 
Since $L(m,n)$ is greater than or equal to both $L(m-1, n)$ and $L(m, n-1)$,  
we see that if $(m,n)\in G_t$, then both points $(m-1, n)$ and $(m,n-1)$ are in $G_t$. 
The set $G_t$ grows with time $t$. 
If we regard $G_t$ in the first quadrant as a stack of boxes, we can add a new box only at the corners. 

A special case is the exponential DLPP in which $w_s$ are exponentially distributed with mean $1$. 
In this case, each corner of $G_t$ grows independently of rate $1$, i.e., a unit box can be added to each corner independent at rate $1$.  
Thus, the boundary of $G_t$ is a rotation of the height function of the corner growth process. 
More precisely, for the TASEP with the step initial condition and the exponential DLPP, 
\beq \label{eq:heightandDLPP}
	\kheight (m-n, t )\ge m+n \quad \text{ if and only if } \quad L(m,n)\le t .
\eeq
If the TASEP starts with a different initial condition, we need to consider the exponential DLPP on a subset of $\Z^2$, determined by the initial condition.

\subsection{Hydrodynamic limit and KPZ limit}

\begin{figure}
\centering
\begin{minipage}{.4\textwidth}
\centering
	\centering\includegraphics[scale=0.14]{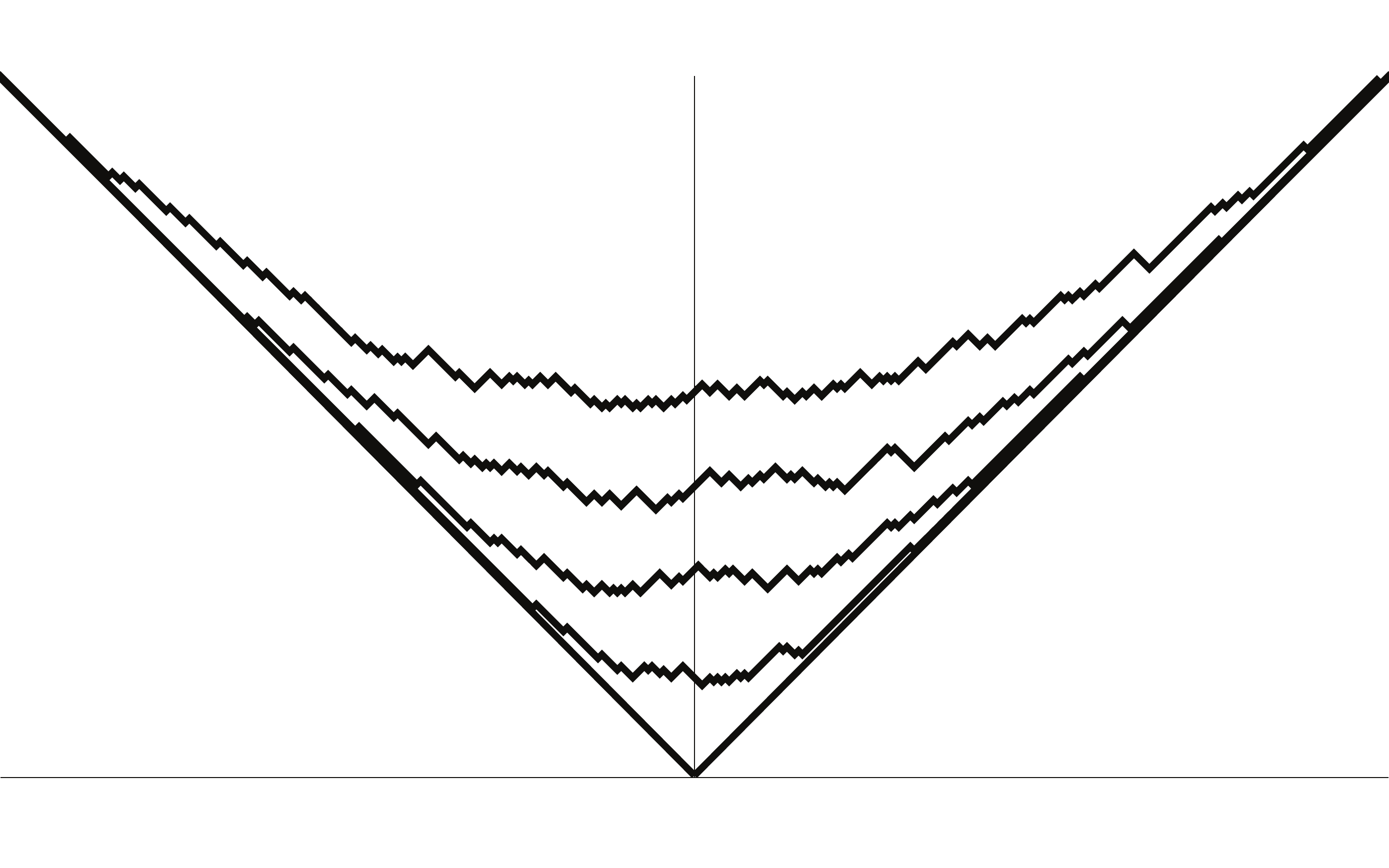}  
\caption{Simulation of the corner growth model}
\label{fig:cornersim}
\end{minipage}
\qquad 
\begin{minipage}{.4\textwidth}
\centering
\begin{tikzpicture}[scale=0.45]
\draw [line width=0.2mm,lightgray, ->] (-4.3,0)--(4.3,0);
\draw [line width=0.2mm,lightgray, ->] (0,0)--(0,4.3);
\node at (4.9, 0) {$x$};
\node at (0.7, 4.3) {$y$};
\draw[domain=-1:1,smooth,variable=\x, thick]  plot ({\x} , {(1+\x*\x)/(2*1)});
\draw[domain=-2:2,smooth,variable=\x, thick]  plot ({\x} , {(2^2+\x*\x)/(2*2)});
\draw[domain=-3:3,smooth,variable=\x, thick]  plot ({\x} , {(3^2+\x*\x)/(2*3)});
\draw[domain=-4:4,smooth,variable=\x, thick]  plot ({\x} , {(4^2+\x*\x)/(2*4)});
\draw[domain= 0:4.4,smooth,variable=\x, thick]  plot ({\x} , {\x});
\draw[domain= -4.4:0,smooth,variable=\x, thick]  plot ({\x} , {-\x});
\end{tikzpicture}
\caption{$y= \bar{\kheight}(x, t)$ for a few values of $t$} 
\label{fig:hbar}
\end{minipage}
\end{figure}

The hydrodynamic limit of TASEP is about $\kheight(x,t)$ when $x$ and $t$ are proportional. 
For the step initial condition, $\kheight(x,0)=|x|$, Rost \cite{Rost81} showed in 1981 that 
\beqq
	\frac{\kheight(xT, tT)}{T}  \to \bar{\kheight}(x,t) 
\eeqq
almost surely as $T\to \infty$, where $\bar{\kheight}(x,t) = \frac{t^2+x^2}{2t}$ for $|x|\le t$ and $\bar{\kheight}(x,t) = |x|$ for $|x|\ge t$. 
See Figure \ref{fig:hbar} for the graph. 
The hydrodynamic limit $\bar{\kheight}$ is deterministic, and it solves Burger's equation \cite{Liggett85, Liggett99}.

The KPZ limit is about the next term, $\kheight(xT, tT)- \bar{\kheight}(x,t)T$. 
Setting $x=0$ for convenience and following the 1:2:3 scale, the KPZ universality conjecture suggests that 
\beqq
	\frac{\kheight(\gamma T^{2/3}, \tau T)- \bar{\kheight}(0,\tau)T}{T^{1/3}}
\eeqq  
converges to a 2d random field. 
If we do not set $x=0$, then we should consider $\kheight(xT+ \gamma T^{2/3}, \tau T)- \bar{\kheight}(x,\tau)T$ in the numerator. 
The limiting 2d field, the KPZ fixed point, depends on the initial condition. The step initial condition for the TASEP becomes the so-called narrow wedge initial condition for the KPZ fixed point. 

TASEP has interpretations as an interacting particle system, a random growth process, and a last passage percolation model. 
Each of these interpretations has natural extensions and generalizations. 
The KPZ universality conjecture is that a large class of models in these generalizations has a universal limit.  
The exact class is not known, but for random growth models, three key features seem to be the locality of growth, some smoothing mechanism, and lateral growths. 
For directed last passage percolation, the universality is expected for all random variables $w_s$ with enough moments and without a large atom at the top of the support of the distribution. The last condition is to prevent the situation that there is always a path connecting $(1,1)$ and $(m,n)$ using only the top value, making the last passage time too concentrated. 

\section{One-point distribution} \label{sec:one}

We discuss one-point KPZ limit theorems from \cite{Baik-Deift-Johansson99} and \cite{Johansson00} mentioned in the introduction, and extensions to other models. 

\subsection{Poisson DLPP}

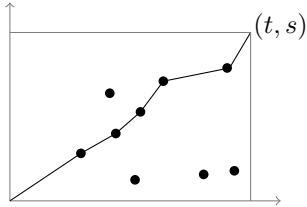
\begin{figure}
\centering \begin{tikzpicture}[scale=0.8]
\draw [->, thin, gray] (0,0) to (4.5,0);
\draw [->, thin, gray] (0,0) to (0, 3.3);
\fill (1.18, 0.79) circle(0.08);
\fill (1.66, 1.79) circle(0.08);
\fill (1.76, 1.12) circle(0.08);
\fill (2.08, 0.35) circle(0.08);
\fill (2.17, 1.48) circle(0.08);
\fill (2.55, 1.99) circle(0.08);
\fill (3.22, 0.44) circle(0.08);
\fill (3.61, 2.21) circle(0.08);
\fill (3.73, 0.50) circle(0.08);
\draw [-, thin] (0,0) to (1.18, 0.79) to (1.76, 1.12) to (2.17, 1.48) to (2.55, 1.99) to (3.65, 2.21) to (4,2.8);
\draw [-, thin, gray] (4,0) to (4, 2.8) to (0, 2.8);
\node at (4.5, 2.9)  {$(t,s)$};
\end{tikzpicture}
\caption{Poisson DLPP}
\label{fig:PoissonDLPP}
\end{figure}

Poisson directed last passage percolation is a variation of the exponential DLPP. 
Consider a realization of a 2d Poisson process in $\R_+^2$. 
An up/right path $p$ this time is defined as the graph of a continuous piecewise linear function of positive slopes connecting Poisson points, as shown in Figure \ref{fig:PoissonDLPP}. 
Let $\pE(p)$ denote the number of the Poisson points on $p$. 
For $(t, s)\in \R_+^2$, define
\beqq
	\pL(t,s) = \sup_{p} \pE(p), 
\eeqq
where the supremum is taken over all up/right paths $p$ from $(0,0)$ to $(t,s)$.
The next theorem follows from \cite{Baik-Deift-Johansson99}. 
The main theorem of \cite{Baik-Deift-Johansson99} is stated for the case of a fixed number of points, 
but the paper proves the Poisson points case first, from which the main theorem follows. 
Since $\pL(t,s) \stackrel{d}{=} \pL(\sqrt{ts}, \sqrt{ts})$, the next result applies to general points $(t,s)$. 


\begin{theorem}[\cite{Baik-Deift-Johansson99}] \label{thm:Poisson}
For every $x\in \R$, 
\beqq
	\lim_{t\to\infty} \Prob\left( \frac{\pL(t,t)-2t}{t^{1/3}}\le x \right) = \FTW (x)
\eeqq
where $\FTW$ is the Tracy-Widom distribution. 
\end{theorem}

The $1/3$-power in $t^{1/3}$ is consistent with the height scale of the KPZ universality: see \eqref{eq:heightandDLPP}. 
This result shows that the one-point marginal of the KPZ fixed point (for the narrow wedge initial condition) must be distributed as the Tracy-Widom distribution. 
The Tracy-Widom distribution is the limiting distribution of the largest eigenvalue of random Hermitian matrices such as Gaussian unitary ensemble matrices \cite{Tracy-Widom94}. 
The connection of the KPZ fixed point and random matrix theory was surprising and unexpected. See Section \ref{sec:TASEPRM} for more on this connection. 


\subsection{Longest increasing subsequence}

The Poisson DLPP is particularly interesting due to its connection to longest increasing subsequences of random permutations. 
Note that finitely many points in a rectangle with distinct $x$ and $y$ coordinates can be associated with a permutation by considering the relative orderings of the coordinates. 
For example, the points in Figure \ref{fig:PoissonDLPP} are associated with the permutation $\pi=475168293$. 
For this permutation, the subsequence $45689$ is an increasing subsequence.  
Furthermore, it is the longest increasing subsequence, and its length, $5$, is equal to the last passage time $L(t,s)$. 

Let $\ell_N$ denote the length of longest increasing subsequences of a uniformly random permutation of size $N$. 
Then, $L(t,s)$ has the same distribution as $\ell_N$, where $N$ is a Poisson random variable of mean $ts$. 
Using this connection, Theorem \ref{thm:Poisson} implies, after a de-Poissonization argument, that 
$\frac{\ell_N- 2\sqrt{N}}{N^{1/6}}$ converges in distribution to the Tracy--Widom distribution.

The problem of determining the large $N$ behavior of $\ell_N$ has a long history. The existence of the almost sure limit of $\ell_N/\sqrt{N}$ was proved by Hammersley in \cite{Hammersley72} using Kingman's ergodic subadditivity theorem. 
The fact that the limit is $2$, known as Ulam's problem, 
was proved independently by two famous papers of Ver\v{s}ik --Kerov \cite{Vershik-Kerov77} and Logan--Shepp \cite{Logan-Shepp77} in 1977. 
However, the limiting distribution and the variance (which is of order $N^{2/3}$) remained an open problem until the work \cite{Baik-Deift-Johansson99}. 
Interested readers are encouraged to consult \cite{Aldous-Diaconis99, Stanley07, Romik15, Baik-Deift-Suidan}.  


\subsection{Exponential DLPP}

Soon after Theorem \ref{thm:Poisson} was proved, Johansson showed that the exponential DLPP model also satisfies a similar limit theorem \cite{Johansson00}. 
We state the result in terms of the height function of the TASEP. 

\begin{theorem}[\cite{Johansson00}] \label{thm:Johansson}
Assume the step initial condition for the TASEP. 
Then, for every $( \kt, \kga, \kh) \in \R_+\times \R \times \R$, 
\beqq
	\lim_{T\to \infty} \Prob \left(  \frac{\kheight(\kga T^{2/3}, 2 \kt T) - \tau T}{-T^{1/3}} \le \kh \right) 
	=  \FTW \left( \frac{\kh}{\kt^{1/3}} + \frac{\kga^2}{4\kt^{4/3}} \right) . 
\eeqq
\end{theorem}

Note from either side of the equation that the limit remains unchanged if we rescale 
\beq \label{eq:rescale}
	(\kh, \kga, \kt ) \mapsto (\alpha \kh, \alpha^2 \kga, \alpha^3 \kt)
\eeq 
for any $\alpha>0$. 




\subsection{Integrable method} \label{sec:integrablemethod}

The above theorems were obtained by explicitly computing the finite-time distribution function and then taking the large limit of the formula. 
In particular, for the TASEP, the finite-time formula is given by the Fredholm determinant of an operator. 
After suitable scaling and conjugation, the operator converges to the so-called Airy operator, the Fredholm determinant of which is the Tracy-Widom distribution. 


Johansson obtained the finite-time distribution formula for TASEP using a combinatorial interpretation similar to the longest increasing subsequence problem and connecting to the so-called Schur measure \cite{Okounkov01}. 
The Schur measure on integer partitions is defined in terms of the Schur function and contains many parameters. 
The one-point distribution of the TASEP arises by taking a special limit of the parameters.   

A different proof computes the transition probabilities of the TASEP explicitly and then takes an appropriate sum over the configuration space to obtain the finite-time distribution.
To find the transition probabilities, we solve the Kolmogorov forward equation, which is a linear differential equation with non-constant coefficients due to the exclusion property of the particles. 
This equation was solved explicitly in \cite{Schutz97} by applying the coordinate Bethe ansatz method from mathematical physics \cite{Gaudin14, Sutherland04}, which consists of changing the Kolmogorov equation to a linear differential equation with constant coefficients (the free evolution equation) but with complicated boundary conditions. 
Taking the sum of transition probabilities over particular configurations is more technical, and this part was done in \cite{Nagao-Sasamoto04, Rakos-Schutz05} to re-derive the result of Johansson. 

Both methods, which are algebraic and exact, are significantly extended to prove a one-point KPZ limit theorem for many other models. 
The following is the list of some of such integrable (exactly solvable) models. Of course, the list and references are far from exhaustive. 

\begin{itemize}
\item Interacting particle systems: PushASEP, ASEP, $q$-TASEP, $q$-Hahn ASEP \cite{Borodin-Ferrari08, Tracy-Widom09, Borodin-Ferrari14, Barraquand-Corwin16}. 
\item Random growth models: KPZ equation, stochastic heat equation  \cite{Amir-Corwin-Quastel11, Borodin-Corwin-Ferrari-Veto15}.
\item DLPP and directed polymers: O'Connell--Yor semi-discrete polymer, log Gamma polymer \cite{O'Connell12, Borodin-Corwin-Remenik13, Borodin-Corwin-Ferrari14}. 
\end{itemize}

The underlying algebraic structures of these integrable models are  generalized greatly by Macdonald processes \cite{Borodin-Corwin13} and the stochastic six vertex model  \cite{Borodin-Corwin-Gorin16, Corwin-Petrov16}. 
They are umbrella models with many parameters whose specializations produce the above models.  Though many, these integrable models are still isolated examples. For instance, DLPP with general random variables, other than exponential and geometric random variables, does not seem to be integrable.  
See Section \ref{sec:univ} for some comments for non-integrable models.


\section{Multi-point distributions}\label{sec:multitime}

Pr\"ahofer and Spohn \cite{Prahofer-Spohn02}, and Johansson \cite{Johansson03} extended the one-point distribution results of Theorem \ref{thm:Poisson} and \ref{thm:Johansson} to equal-time, multi-position distributions for  the Poisson DLPP and the TASEP with step initial condition, respectively. 
Their results were further extended to other models and initial conditions by \cite{Sasamoto05, Borodin-Ferrari-Prahofer-Sasamoto07, Borodin-Ferrari-Sasamoto08, Borodin-Ferrari-Sasamoto08a} in 2005-2008.  
These results confirmed, in particular, that the spatial correlations are of order $T^{2/3}$ and identified the equal-time slice of the KPZ fixed point for several initial conditions. 
See also \cite{Dimitrov20} for more recent progress on other more difficult models. 



On the other hand, multi-time distributions and fully 2d multi-position distributions remained unanswered for a while
though some short and long time correlations were studied in \cite{Ferrari-Spohn16}, confirming that the time scale is $T$.    
In a breakthrough paper \cite{Matetski-Quastel-Remenik17} published in 2021, Mateski, Quastel, and Remenik proved the convergence of the entire 2d height field of the TASEP. 
The limiting 2d field, the KPZ fixed point, is constructed as a Markov process with explicit transition probabilities. 
The result applies to general initial conditions. 
The authors used the result of \cite{Sasamoto05, Borodin-Ferrari-Prahofer-Sasamoto07} on the transition probabilities of the TASEP for general initial conditions and proved that they converge. 
Dauvergne, Ortmann, and Vir\'ag gave an alternative formulation of the KPZ fixed point in terms of a variational formula and proved the field convergence for another model, the Brownian DLPP \cite{Dauvergne-Ortmann-Virag18}.  
See also \cite{Virag20}.

In the meantime, Johansson and Rahman \cite{Johansson-Rahman19}, and Liu \cite{Liu19} computed the limit of 2d multi-point distributions of the discrete-time PNG and the continuous-time TASEP, respectively, in 2021 and 2022. 
Their results give an explicit formula of multi-point distributions for the KPZ fixed point with the narrow wedge initial condition.  
See Section \ref{sec:dfformula} for the formula. 
Two-time distributions were previously computed in  \cite{Johansson17, Johansson19}. 

\section{Half-infinite space} \label{sec:half}

We discussed so far models on infinite spaces. 
For example, the TASEP was defined on $\Z$. 
In this and the following sections, we consider different domains and their effects on the limit. 


Consider the TASEP on the half-infinite space $\Z_+\cup\{0\}$. 
We introduce a parameter $\alpha>0$ representing the injection rate at site $0$: if the origin is empty, a new particle is injected with rate $\alpha$. Once injected, particles follow the usual TASEP rule. 
Suppose that we start with the empty configuration. 
If we could inject particles freely without being blocked by existing particles in the domain, then the height function at the origin would satisfy $\kheight(0,T)/T\to 2\alpha$ in probability as $T\to \infty$. 
However, due to the particles already in the domain, the height grows at a slower rate, and the hydrodynamic limit at the origin turns out to be 
\beqq
	\frac{\kheight (0, T)}{T} \to \max \big\{ 2\alpha(1-\alpha), \frac{1}{2} \big\}, 
\eeqq
showing that the effective injection rate is $\max\{ \alpha(1-\alpha), 1/4\}$. 
The formula changes at $\alpha=1/2$. 

The papers \cite{Baik-Rains01a, Baik-Rains01} obtained a one-point KPZ limit theorem at the origin for the Poisson PNG and discrete-time TASEP models. 
The height scales as $T^{1/3}$ for $\alpha\ge 1/2$ and as $T^{1/2}$ for $\alpha<1/2$. 
The limiting distribution is  a variation of the Tracy-Widom distribution for $\alpha>1/2$, another variation for $\alpha=1/2$, and the Gaussian distribution for $\alpha<1/2$. 
The result is extended to general positions and equal-time, multi-position distributions in \cite{Sasamoto-Imamura04, Baik-Barraquand-Corwin-Suidan18} for the Poisson PNG and the discrete and continuous-time TASEP. 
The one-point distribution result was also extended to other models recently. The ASEP was analyzed in \cite{Barraquand-Borodin-Corwin-Wheeler18} for a particular value of the boundary condition.
The paper \cite{Imamura-Mucciconi-Sasamoto22} established one-point distribution results for log-Gamma directed polymer, q-PushTASEP, and multiplicative stochastic heat equation,
 for general values of the boundary condition.  However, multi-time limit theorems are still yet to be established

\section{Ring domain} \label{sec:ring}

Consider the TASEP on the integer ring $\Z_L=\Z/L\Z=\{0, 1, \cdots, L-1\}$ where we identify sites $L$ and $0$. 
An equivalent model is the TASEP on $\Z$ that is spatially periodic, which we may call the periodic TASEP. 
We call $L$ the size of the ring or the period of the periodic TASEP.

\begin{figure}\centering
	\begin{minipage}{.45\textwidth} \centering
\begin{tikzpicture}[scale=0.3]
\draw (0, 4) circle(0.2);
\draw (2, 3.464) circle(0.2);
\draw (3.464, 2) circle(0.2);
\draw (4, 0) circle(0.2);
\draw (3.464, -2) circle(0.2);
\draw (2, -3.464) circle(0.2);
\fill (0, -4) circle(0.2);
\fill (-2, -3.464) circle(0.2);
\fill (-3.464, -2) circle(0.2);
\fill (-4, 0) circle(0.2);
\fill (-3.464, 2) circle(0.2);
\fill (-2, 3.464) circle(0.2);
\end{tikzpicture}
\caption{Step initial condition on an integer ring}
\label{fig:Ring}
	\end{minipage}
	\begin{minipage}{.45\textwidth} \centering
\begin{tikzpicture}[scale=1.2]
\draw [line width=0.2mm,lightgray] (-1.6,0)--(1.6,0);
\draw [line width=0.2mm,lightgray] (0, 0)--(0,1.7);
\draw[domain= -1.5:-1,smooth,variable=\x]  plot ({\x} , {-\x-1});
\draw[domain= -1:-0.5,smooth,variable=\x]  plot ({\x} , {\x+1});
\draw[domain= -0.5:0,smooth,variable=\x]  plot ({\x} , {-\x});
\draw[domain= 0:0.5,smooth,variable=\x]  plot ({\x} , {\x});
\draw[domain= 0.5:1,smooth,variable=\x]  plot ({\x} , {-\x+1});
\draw[domain= 1:1.5,smooth,variable=\x]  plot ({\x} , {\x-1});
\draw[domain=-1.5:-0.5,smooth,variable=\x ]  plot ({\x} , {(0.5^2+(\x+1)*(\x+1))/(2*1/2)});
\draw[domain=-0.5:0.5,smooth,variable=\x ]  plot ({\x} , {(0.5^2+\x*\x)/(2*1/2)});
\draw[domain=0.5:1.5,smooth,variable=\x ]  plot ({\x} , {(0.5^2+(\x-1)*(\x-1))/(2*1/2)});
\draw[domain=-1.5:-0.5,smooth,variable=\x ]  plot ({\x} , {(1^2+(\x+1)*(\x+1))/(2*1)});
\draw[domain=-0.5:0.5,smooth,variable=\x ]  plot ({\x} , {(1^2+\x*\x)/(2*1)});
\draw[domain=0.5:1.5,smooth,variable=\x ]  plot ({\x} , {(1^2+(\x-1)*(\x-1))/(2*1)});
\draw[domain=-1.5:-0.5,smooth,variable=\x ]  plot ({\x} , {(1.5^2+(\x+1)*(\x+1))/(2*1.5)});
\draw[domain=-0.5:0.5,smooth,variable=\x ]  plot ({\x} , {(1.5^2+\x*\x)/(2*1.5)});
\draw[domain=0.5:1.5,smooth,variable=\x ]  plot ({\x} , {(1.5^2+(\x-1)*(\x-1))/(2*1.5)});
\draw[domain=-1.5:-0.5,smooth,variable=\x ]  plot ({\x} , {(2^2+(\x+1)*(\x+1))/(2*2)});
\draw[domain=-0.5:0.5,smooth,variable=\x ]  plot ({\x} , {(2^2+\x*\x)/(2*2)});
\draw[domain=0.5:1.5,smooth,variable=\x ]  plot ({\x} , {(2^2+(\x-1)*(\x-1))/(2*2)});
\draw[domain=-1.5:-0.5,smooth,variable=\x ]  plot ({\x} , {(2.5^2+(\x+1)*(\x+1))/(2*2.5)});
\draw[domain=-0.5:0.5,smooth,variable=\x ]  plot ({\x} , {(2.5^2+\x*\x)/(2*2.5)});
\draw[domain=0.5:1.5,smooth,variable=\x ]  plot ({\x} , {(2.5^2+(\x-1)*(\x-1))/(2*2.5)});
\draw[domain=-1.5:-0.5,smooth,variable=\x ]  plot ({\x} , {(3^2+(\x+1)*(\x+1))/(2*3)});
\draw[domain=-0.5:0.5,smooth,variable=\x ]  plot ({\x} , {(3^2+\x*\x)/(2*3)});
\draw[domain=0.5:1.5,smooth,variable=\x ]  plot ({\x} , {(3^2+(\x-1)*(\x-1))/(2*3)});
\end{tikzpicture}
\caption{Hydrodynamic limits of periodic height function when $t=O(t)$}
\label{fig:Ringmaro}
	\end{minipage}
\end{figure}

The number $N$ of particles in the TASEP on the ring is preserved. 
We assume the step initial condition shown in Figure \ref{fig:Ring}. 
For the convenience of presentation, we assume that $L$ is an even integer and $N=L/2$ so that the particle density $\rho=N/L=1/2$. 
The initial height function of the periodic TASEP is the bottom curve in Figure \ref{fig:Ringmaro}. 
The other curves are the hydrodynamic limits, as time $t$ and period $L$ tend to infinity proportionally when $t/L=0.5n$ for $n=1,2, \cdots, 6$. 

Consider two cases, one that $t\to\infty$ with $L$ fixed and the other that $L\to\infty$ with $t$ fixed. 
If $t\to\infty$ with $L$ fixed, the periodic corner growth model becomes essentially a one-dimensional growth model, and we expect that the height scales as $t^{1/2}$ and converges to the Gaussian distribution. 
On the other hand, if $L\to \infty$ with $t$ fixed, then the periodic TASEP becomes the usual TASEP on $\Z$. 
Thus, the height scales as $t^{1/3}$ if we let $t\to\infty$ after taking $L\to \infty$. 
An interesting intermediate regime is when $L, t\to\infty$ simultaneously such that 
\beq \label{eq:relax}
	t=O(L^{3/2}). 
\eeq
Since the spatial scale for KPZ limit theorems on infinite spaces is $t^{2/3}$, we expect that the height functions at all positions on the ring of size $L$ are correlated non-trivially. 
If \eqref{eq:relax} holds, we say that we are in the relaxation time regime.

There were some results on transition probabilities and the spectral gap of the generator in the relaxation time regime, such  as \cite{Gwa-SpohnBethe}.  
However, KPZ limit theorems 
were obtained only recently. 
The paper \cite{Prolhac16}, which is a non-completely rigorous physics paper, and \cite{Baik-Liu18} obtained a one-point KPZ limit theorem almost at the same time independently. 
This result was further extended to 2d multi-point distributions in \cite{Baik-Liu19, Baik-Liu21}. 

\begin{theorem}[\cite{Baik-Liu19}] \label{thm:PTASEP}
Consider the TASEP on a ring of size $L$ with the step initial condition and extend it to the periodic TASEP. 
Assume that $L$ is even and $\rho=N/L= 1/2$ for the convenience of presentation. 
Set 
\beqq
	T=L^{3/2}. 
\eeqq
For $i=1, \cdots, m$, let 
 $(\pkga_i,  \pkt_i, \pkh_i) \in \R \times  \R_+\times \R$ and assume that $\pkt_1<\cdots< \pkt_m$. 
Then, 
\beqq
	\lim_{T\to \infty} \Prob \Big( \bigcap_{i=1}^m  \big\{  \frac{\pkheight(\pkga_i T^{2/3}, 2 \pkt_i T) - \pkt_i T}{-T^{1/3}} \le \pkh_i   \big\} \Big) 
	=   \pkFdist_m(\pkh; \pkga, \pkt)
\eeqq
for an $m$-point distribution function $\pkFdist_m$ described in the next section. 
\end{theorem}

Like the infinite space case, we expect that the height field of the periodic models in the relaxation time regime converges to a universal field, which we may call the periodic KPZ fixed point. 
The function $\pkFdist_m$ should be the $m$-point distribution of this conjectured periodic KPZ fixed point with the (periodic) narrow wedge initial condition. 
It is naturally periodic with respect to $\pkga_i\mapsto \pkga_i+1$. 
However, unlike the KPZ fixed point, $\pkFdist_m$ is not invariant under the rescaling \eqref{eq:rescale}. 
Indeed, we conjecture that $\pkFdist_m$ interpolates the KPZ fixed point and one-dimensional Brownian motion. Concretely, we expect that  
\beqq
	 \lim_{\epsilon\to 0} \pkFdist_m(\pkh^\epsilon; \pkga^\epsilon, \pkt^\epsilon) =  \kFdist_m(\pkh; \pkga, \pkt) 
	 \qquad \text{where $(\pkh_i^\epsilon, \pkga_i^\epsilon, \pkt_i^\epsilon) = ( (\pkt_i \epsilon)^{1/3}\pkh_i, (\pkt_i \epsilon)^{2/3} \pkga_i,  \pkt_i \epsilon)$} 
\eeqq
and $\kFdist_m$ is the $m$-point distribution of the KPZ fixed point, and that 
\beqq
	\lim_{s\to \infty} \pkFdist_m ( \pkh^s; \pkga,  \pkt^s  )  = G_m(\pkh; \pkt) 
	\qquad \text{where $( \pkh_i^s, \pkt^s) = ( - s \pkt_i + \frac{s^{1/2} \pi^{1/4} }{\sqrt{2}} \pkh_i, s \pkt_i)$ }
\eeqq
and $G_m$ is the $m$-point distribution of a Brownian motion at times $\pkt_1, \cdots, \pkt_m$. 
These conjectures were proved for $m=1$ in \cite{Baik-Liu-Silva20}, assuming $\pkga_1=0$ for the $\epsilon\to 0$ case.



Theorem \ref{thm:PTASEP} is also proved for the discrete-time TASEP on a ring \cite{Liao21}. 
However,  extending the result to other integrable models is at the moment a challenge.

\section{Formula of multi-point distribution functions} \label{sec:dfformula}

\subsection{Formula for KPZ fixed point}

Let $\kFdist_m$ be the $m$-point distribution of the KPZ fixed point with the narrow wedge initial condition. The result of  \cite{Liu19} implies that  
\beq \label{eq:KPZformula}
	\kFdist_m(\kh; \kga, \kt) 
	= \frac1{(2\pi \ii)^{m-1}} \oint \cdots \oint  \frac{\det( \kid- \kK_\zeta)}{\zeta_1(1-\zeta_1) \cdots \zeta_{m-1} (1-\zeta_{m-1})} \dd \zeta_1 \cdots \dd \zeta_{m-1}
\eeq
where $\zeta=(\zeta_1, \cdots, \zeta_{m-1})$, and the contours are nested circles of radii less than $1$ centered at the origin. 
The operator $\kK_\zeta$ acts on $L^2(\Sigma)$, where 
$\Sigma$ is the union of $4m-2$ contours in Figure \ref{fig:KPZcontour} that extend to infinity with angle $\pi/5$ from the $x$-axis.
The kernel of $\kK_\zeta$ can be written \cite{Baik-Prokhorov-Silva21}  as a simple conjugation of the kernel 
\beq \label{eq:kpzKkernel}
	\kK^{\mathrm{conj}}_\zeta(\ku, \kv) = \frac{\ka(\ku)^T \kD(\ku)^T \kD(\kv) \kb( \kv)}{\ku-\kv}, \qquad \ku, \kv\in \Sigma, 
\eeq
which is zero for $\ku=\kv$. 
The $(m+1)\times (m+1)$ matrix 
\beqq
	\kD(z)= \diag\big( e^{- \frac{1}{3} \kt_1 z^3+ \frac{1}{2} \kga_1 z^2 + \kh_1 z}, \cdots, e^{- \frac{1}{3} \kt_m z^3 + \frac{1}{2} \kga_m z^2 + \kh_m z}, 1 \big) . 
\eeqq
The $(m+1)\times 1$ vectors $\ka(z)$ and $\kb(z)$ are simple and explicit, and they do not depend on $\kh_i, \kga_i, \kt_i$. 
Note that the exponent $- \frac{1}{3} \kt_i z^3+ \frac{1}{2} \kga_i z^2+ \kh_i z$ in $\kD(z)$ is unchanged if we rescale as \eqref{eq:rescale} and $z\mapsto \alpha^{-1} z$. 
This is consistent with the rescaling property of the KPZ fixed point. 

\begin{figure}
\centering
\begin{minipage}{.4\textwidth}\centering
\begin{tikzpicture}[scale=0.3]
\draw [line width=0.2mm,lightgray] (-7.3,0)--(7.3,0);
\draw [line width=0.2mm,lightgray] (0,-4.3)--(0,4.3);
\draw[domain=-3:3,smooth,variable=\y, thick]  plot ({-(2*\y*\y+1)^(1/2)},{\y});
\draw[domain=-3:3,smooth,variable=\y, thick]  plot ({-1-(2*\y*\y+1)^(1/2)},{\y});
\draw[domain=-3:3,smooth,variable=\y, thick]  plot ({-2-(2*\y*\y+1)^(1/2)},{\y});
\draw[domain=-3:3,smooth,variable=\y, thick]  plot ({-3-(2*\y*\y+1)^(1/2)},{\y});
\draw[domain=-3:3,smooth,variable=\y, thick]  plot ({-4-(2*\y*\y+1)^(1/2)},{\y});
\draw[domain=-3:3,smooth,variable=\y, thick]   plot ({(2*\y*\y+1)^(1/2)},{\y});
\draw[domain=-3:3,smooth,variable=\y, thick]  plot ({1+(2*\y*\y+1)^(1/2)},{\y});
\draw[domain=-3:3,smooth,variable=\y, thick]   plot ({2+(2*\y*\y+1)^(1/2)},{\y});
\draw[domain=-3:3,smooth,variable=\y, thick]  plot ({3+(2*\y*\y+1)^(1/2)},{\y});
\draw[domain=-3:3,smooth,variable=\y, thick]   plot ({4+(2*\y*\y+1)^(1/2)},{\y});
\end{tikzpicture}
\caption{The space $\Sigma$ for KPZ fixed point when $m=3$}
\label{fig:KPZcontour}
\end{minipage}
\quad
\begin{minipage}{.45\textwidth}\centering
\begin{tikzpicture}[scale=0.3]
\draw [lightgray] (-5,0)--(5,0);
\draw [lightgray] (0,-4)--(0,4);
\draw[domain=-4:4,smooth,variable=\y,lightgray]  plot ({-(\y*\y+1.02)^(1/2)},{\y});   
\draw[domain=-4:4,smooth,variable=\y,lightgray]  plot ({-(\y*\y+4.61)^(1/2)},{\y});   
\draw[domain=-4:4,smooth,variable=\y,lightgray]  plot ({-(\y*\y+9.21)^(1/2)},{\y});   
\draw[domain=-4:4,smooth,variable=\y,lightgray]  plot ({(\y*\y+1.02)^(1/2)},{\y});   
\draw[domain=-4:4,smooth,variable=\y,lightgray]  plot ({(\y*\y+4.61)^(1/2)},{\y});   
\draw[domain=-4:4,smooth,variable=\y,lightgray]  plot ({(\y*\y+9.21)^(1/2)},{\y});   
\fill (-3.93, -3.30) circle[radius=3.5pt] node [above,shift={(0pt,0pt)}] {};   
\fill (-3.06, -2.18) circle[radius=3.5pt] node [above,shift={(0pt,0pt)}] {};
\fill (-2.15, -0.19) circle[radius=3.5pt] node [above,shift={(0pt,0pt)}] {};
\fill (-2.93, 2.00) circle[radius=3.5pt] node [above,shift={(0pt,0pt)}] {};
\fill (-3.83, 3.17) circle[radius=3.5pt] node [above,shift={(0pt,0pt)}] {};
\fill (-3.67, -3.53) circle[radius=3.5pt] node [above,shift={(0pt,0pt)}] {};   
\fill (-2.67, -2.49) circle[radius=3.5pt] node [above,shift={(0pt,0pt)}] {};
\fill (-1.08, -0.37) circle[radius=3.5pt] node [above,shift={(0pt,0pt)}] {};
\fill (-2.53, 2.32) circle[radius=3.5pt] node [above,shift={(0pt,0pt)}] {};
\fill  (-3.56, 3.42) circle[radius=3.5pt] node [above,shift={(0pt,0pt)}] {};
\fill (-4.29, -3.03) circle[radius=3.5pt] node [above,shift={(0pt,0pt)}] {};    
\fill (-3.57, -1.87) circle[radius=3.5pt] node [above,shift={(0pt,0pt)}] {};
\fill (-3.04, -0.13) circle[radius=3.5pt] node [above,shift={(0pt,0pt)}] {};
\fill (-3.48, 1.69) circle[radius=3.5pt] node [above,shift={(0pt,0pt)}] {};
\fill  (-4.20, 2.90) circle[radius=3.5pt] node [above,shift={(0pt,0pt)}] {};
\fill (3.93, 3.30) circle[radius=3.5pt] node [above,shift={(0pt,0pt)}] {};   
\fill (3.06, 2.18) circle[radius=3.5pt] node [above,shift={(0pt,0pt)}] {};
\fill (2.15, 0.19) circle[radius=3.5pt] node [above,shift={(0pt,0pt)}] {};
\fill (2.93, -2.00) circle[radius=3.5pt] node [above,shift={(0pt,0pt)}] {};
\fill (3.83, -3.17) circle[radius=3.5pt] node [above,shift={(0pt,0pt)}] {};
\fill (3.67, 3.53) circle[radius=3.5pt] node [above,shift={(0pt,0pt)}] {};   
\fill (2.67, 2.49) circle[radius=3.5pt] node [above,shift={(0pt,0pt)}] {};
\fill (1.08, 0.37) circle[radius=3.5pt] node [above,shift={(0pt,0pt)}] {};
\fill (2.53, -2.32) circle[radius=3.5pt] node [above,shift={(0pt,0pt)}] {};
\fill  (3.56, -3.42) circle[radius=3.5pt] node [above,shift={(0pt,0pt)}] {};
\fill (4.29, 3.03) circle[radius=3.5pt] node [above,shift={(0pt,0pt)}] {};    
\fill (3.57, 1.87) circle[radius=3.5pt] node [above,shift={(0pt,0pt)}] {};
\fill (3.04, 0.13) circle[radius=3.5pt] node [above,shift={(0pt,0pt)}] {};
\fill (3.48, -1.69) circle[radius=3.5pt] node [above,shift={(0pt,0pt)}] {};
\fill  (4.20, -2.90) circle[radius=3.5pt] node [above,shift={(0pt,0pt)}] {};
\end{tikzpicture}
\caption{The space $\pspace$ for periodic KPZ fixed when $m=3$.}
\label{fig:periodicKPZset}
\end{minipage}
\end{figure}

\subsection{Periodic case}

The formula of \cite{Baik-Liu19} is 
\beq \label{eq:PKPZformula}
	\pkFdist_m(\pkh; \pkga, \pkt)
	= \oint \cdots \oint  \pkCdisc(\zeta) \det(\pid-\pK_\zeta)   \dd \zeta_1 \cdots \dd \zeta_{m}, 
\eeq
where this time there are $m$ integrals and $\pkCdisc(\zeta)$ is an explicit function expressed in terms of polylog functions. 
The kernel of the operator $\pK_\zeta$ is of the same form as \eqref{eq:kpzKkernel} but $\ka(z)$ and $\kb(z)$ are slightly different. 
The key change is the space for $\pK_{\zeta}$. It is $\ell^2(\pspace)$, where
$\pspace= \pspace_1\cup\cdots \cup \pspace_m$ and $\pspace_i$ is the discrete set of the roots of the equation 
\beq \label{eq:Betheeqlimit}
	e^{-s^2/2}= \zeta_i , 
\eeq
shown in Figure \ref{fig:periodicKPZset}. 
See the following subsection for how this equation arises.

\subsection{Transition probabilities} \label{sec:Bethe}

We mentioned in Subsection \ref{sec:integrablemethod} that one way of proving a KPZ limit theorem for the TASEP is to compute the transition probabilities explicitly and then take an appropriate sum to find the finite-time distribution functions. 
The summation part is often more technical, but here we discuss the transition probabilities to see how \eqref{eq:Betheeqlimit} arises. 
 
Suppose that there are only $N$ particles for the TASEP on $\Z$. 
Let $\mathcal{W}_N= \{ (a_1, \cdots, a_N)\in \Z^N: a_1<\cdots< a_N\}$ be the ordered set of the particle locations. 
Sch\"utz \cite{Schutz97} showed that 
\beq \label{eq:tranpro}
	\Prob_{Y}(X; t) =  \det \left[ \frac1{2\pi \ii} \oint s^{j-i+1}(s+1)^{-x_i+y_j+i-j} e^{ts} \dd s \right]_{i,j=1}^N 
\eeq
for $X$ and $Y$ in $\mathcal{W}_N$, where the contour is a circle that encloses the points $s=0,-1$. 

\begin{figure}
\centering
\begin{minipage}{.25\textwidth}
\centering
\includegraphics[width=3cm]{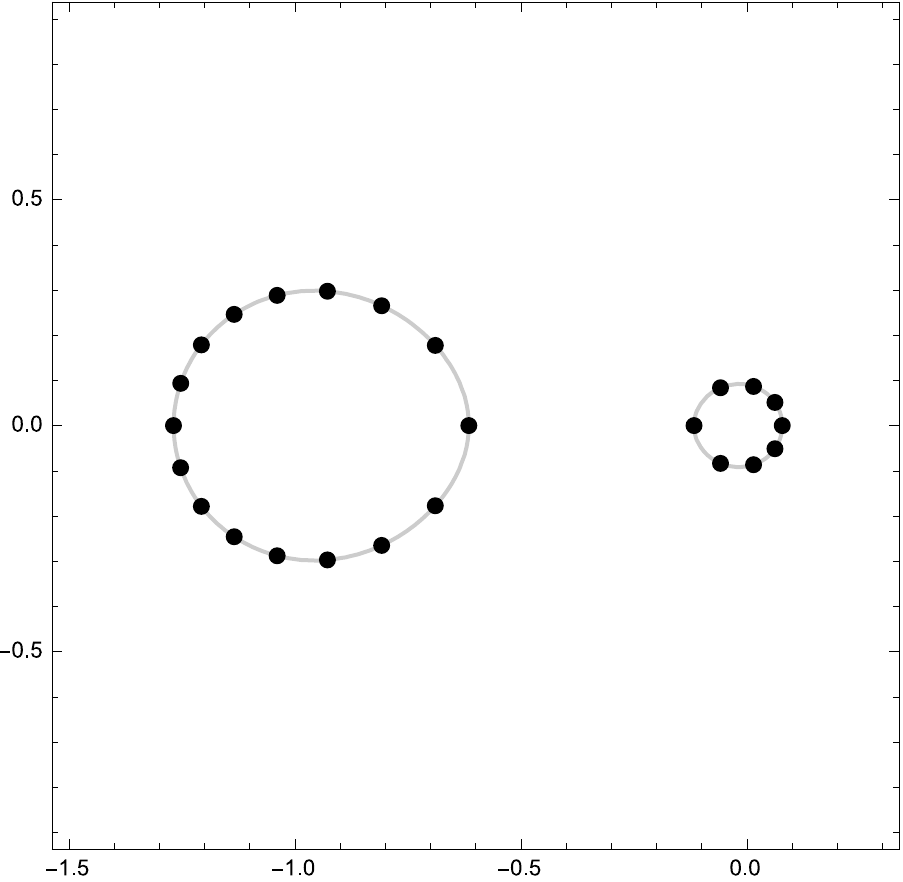}
\end{minipage}
\begin{minipage}{.25\textwidth}
\centering
\includegraphics[width=3cm]{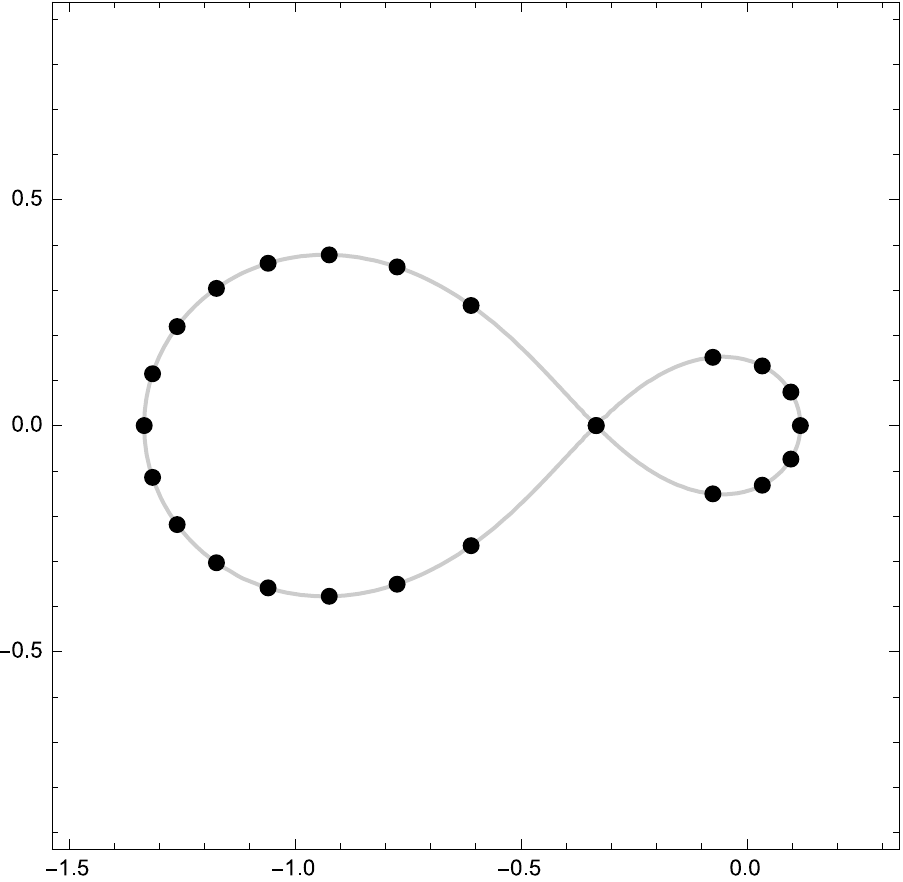}
\end{minipage}
\begin{minipage}{.25\textwidth}
\centering
\includegraphics[width=3cm]{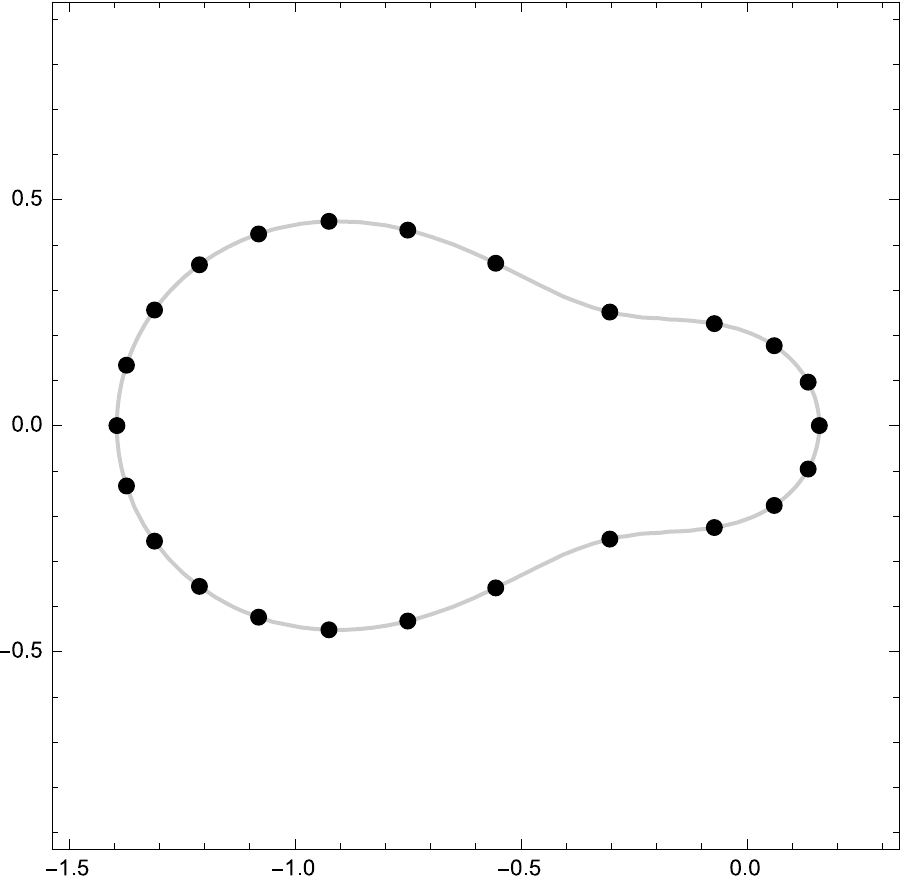}
\end{minipage}
\caption{Bethe roots when $L=24$ and $N=8$ for three values of $z$}
\label{fig:Betheroots}
\end{figure}

For the periodic TASEP, the particle locations can be represented by the set 
\beqq
	\mathcal{W}_N^L=  \{ (a_1, \cdots, a_N)\in \Z^N: a_1<\cdots< a_N <a_1+L\} . 
\eeqq
Note that if we consider the TASEP on a ring, this set keeps track of global circulations of the particles. 
We showed in \cite{Baik-Liu18} that 
\beq \label{eq:transprofper}
	\Prob_{Y}(X; t) = \oint \det \left[ \frac1{L} \sum_{w} \frac{w^{j-i+1}(w+1)^{-x_i+y_j+i-j} e^{tw}}{w+N/L} \right]_{i,j=1}^N  \frac{\dd z}{2\pi \ii z}
\eeq
for $X, Y\in \mathcal{W}_N^L$, where the integral contour for $z$ is any circle enclosing the origin, and the sum inside is over the roots of the equation 
\beq \label{eq:Betheequation}
	w^N (w+1)^{L-N}=z . 
\eeq 
See Figure \ref{fig:Betheroots}. 

We now explain the equation \eqref{eq:Betheequation}. 
Due to the periodicity, if $A=(a_1, \cdots, a_N)$ is in $\mathcal{W}_N(L)$, then $A'=(a_2, \cdots, a_N, a_1+L)$ also represents the same particle configuration of the periodic TASEP. 
Thus, the transition probability should remain the same if we replace $(x_1, \cdots, x_N)$ and $(y_1, \cdots, y_N)$ by $(x_2, \cdots, x_N, x_1+L)$ and $(y_2, \cdots, y_N, y_1+L)$, respectively. 
We can check directly that the determinant in \eqref{eq:transprofper} is unchanged thanks to the equation \eqref{eq:Betheequation}. 
Equation \eqref{eq:Betheequation} takes care of the labeling ambiguity in the periodic case. 
Here the variable $z$ could have been any fixed constant, but making it as a free parameter turns out to be the right choice. 
Since \eqref{eq:tranpro} and \eqref{eq:transprofper} were found by solving the Kolmogorov forward equation using the Bethe ansatz method, the roots of \eqref{eq:Betheequation} are called the Bethe roots. 

Now, if we set $L=2N$, $w= -\frac12 + \frac{s}{2\sqrt{2N}}$, and $z= (-4)^{-N} \zeta$ in \eqref{eq:Betheequation}, and let $N\to \infty$, then the equation becomes $e^{-s^2/2}=\zeta$, which is \eqref{eq:Betheeqlimit}.


\section{Integrable differential equations} \label{sec:integrableDE}

Distribution functions of the KPZ fixed point have connections to deterministic integrable differential equations. 
As proved in 1994 \cite{Tracy-Widom94}, the Tracy-Widom distribution is expressible in terms of the Painlev\'e II equation, one of a family of six special nonlinear ordinary differential equations  \cite{Fokas-Its-Kitaev-Novokshenov06}. 
The papers \cite{Tracy-Widom03, Adler-van_Moerbeke05, Wang08a, Bertola-Cafasso-2012b, Quastel-Remenik19b} also found differential equations for equal-time, multi-position distribution functions of the KPZ fixed point. 
We state the following result for multi-point distributions for both infinite and periodic domains. 


Define the parameters
\beqq
	\knt_i= \kt_i/3, 
	\qquad \kny_i=\kga_i, 
	\qquad \knx_i= \kh_i, 
\eeqq
and let
\beqq
	\partial_{\knt}= \sum_{i=1}^m \partial_{\knt_i}, \qquad
	\partial_{\kny}= \sum_{i=1}^m \partial_{\kny_i}, \qquad
	\partial_{\knx}= \sum_{i=1}^m \partial_{\knx_i}. 
\eeqq

\begin{theorem}[\cite{Baik-Liu-Silva20, Baik-Prokhorov-Silva21}] \label{thm:mKdV}
Let $\kK=\kK_\zeta$ be the operator in either \eqref{eq:KPZformula} or \eqref{eq:PKPZformula}. 
If $\det(\kid- \kK)\neq 0$, which holds for all but at most countably many parameters, then 
\beqq
	\partial^2_{\knx} \log \det(\kid- \kK) = - \kr^T\kp
\eeqq
for complex-valued $m\times 1$ vector functions $\kp(\knt, \kny, \knx)$ and $\kr(\knt, \kny, \knx)$ which satisfy the equations
\beq \label{eq:NLS}
	\partial_{\kny} \kp = \frac12 \partial_{\knx}^2 \kp -  \kp \kr^T \kp, \qquad
	\partial_{\kny} \kr = -\frac12 \partial_{\knx}^2 \kr +  \kr \kp^T \kr
\eeq
and  
\beq \label{eq:mKdV}
	\partial_{\knt} \kp + \partial_{\knx}^3 \kp - 3 (\partial_{\knx} \kp) \kr^T \kp  - 3 \kp \kr^T (\partial_{\knx} \kp) =0 , \quad
	\partial_{\knt} \kr + \partial_{\knx}^3 \kr - 3  (\partial_{\knx} \kr) \kp^T \kr -  3 \kr \kp^T (\partial_{\knx} \kr ) =0. 
\eeq
\end{theorem}

Equation \eqref{eq:mKdV} is a coupled system of vector-valued modified Korteweg-De Vries (mKdV) equations. 
The scalar mKdV equation is 
$\partial_{\knt} f + \partial_{\knx}^3 f - 6 (\partial_{\knx} f) f^2  =0$. 
On the other hand,  equation \eqref{eq:NLS} forms a coupled system of vector-valued nonlinear forward and backward heat equations. 
They become vector-valued nonlinear Sch{\"o}dinger (NLS) equations if we change $\kny_i\mapsto \ii \kny_i$. 
NLS and mKdV equations are two of the most famous integrable partial differential equations \cite{Ablowitz-Clarkson91}. 
The above two systems of equations can be combined to the Kadomtsev-Petviashvili (KP) equation, another integrable differential equation in three variables.
Theorem \ref{thm:mKdV} was obtained using the fact that the operator is a so-called integrable operator \cite{IIKS, Deift99b}. 

Integrable differential equations have a long history starting with the work of Gardner, Greene, Kruskal, and Miura in 1967. 
They found a scattering transform method, an equation-specific nonlinear Fourier transform, to solve the Korteweg-De Vries equation.  
Like the integrable models in the KPZ universality class, integrable differentiable equations are isolated examples of nonlinear differential equations that can often be solved explicitly and analyzed asymptotically. See, for example, \cite{Ablowitz-Clarkson91, Deift19}. 
It is intriguing that integrable probability models are related to integrable differential equations.


\section{Comments on universality} \label{sec:univ}

KPZ limit theorems are proved for many isolated examples of integrable models. 
In this final section, we discuss a few instances that universality is proved.

\subsection{Thin DLPP}

The universality should hold for the last passage time $L(n,k)$ as $n,k\to \infty$.
It is easy to prove it for thin rectangles. 
Recall from Subsection \ref{sec:DLPP} that $w_s$ denotes the random variable at site $s\in \Z$, representing the passage time through the site. 

\begin{theorem}[\cite{Baik-Suidan05, Bodineau-Martin05}]  \label{thm:thin}
Suppose that $w_s$ is an arbitrary random variable which has all moments. Assume that the mean is zero and the variance is one. Then, for every $x$, 
\beqq
	\lim_{n, k \to \infty} \Prob\Big( \frac{L(n, k )- 2 \sqrt{nk}}{n^{\frac12} k^{-\frac16}} \le x \Big) = \FTW(x) \quad \text{for $k=[n^a]$}
\eeqq
with any $0<a<3/7$. 
\end{theorem}

The restriction $a<3/7$ is technical. If we assume only finite $p$ moments for $p>2$, then the result holds for thinner rectangles satisfying $0<a<\frac{3(p-2)}{7 p}$. 

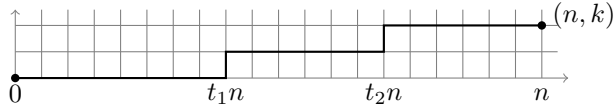
\begin{figure}
\centering
\begin{tikzpicture}[scale=0.7]
\draw [->, thin][gray] (0,0) to (10.5,0);
\draw [->, thin][gray] (0,0) to (0, 1.3);
\draw [-, thin, gray] (0,0.5) to (10.3,0.5);
\draw [-, thin, gray] (0,1) to (10.3,1);
\draw [-, thin, gray] (0.5,0) to (0.5,1.3);
\draw [-, thin, gray] (1,0) to (1,1.3);
\draw [-, thin, gray] (1.5,0) to (1.5,1.3);
\draw [-, thin, gray] (2,0) to (2,1.3);
\draw [-, thin, gray] (2.5,0) to (2.5,1.3);
\draw [-, thin, gray] (3,0) to (3,1.3);
\draw [-, thin, gray] (3.5,0) to (3.5,1.3);
\draw [-, thin, gray] (4,0) to (4,1.3);
\draw [-, thin, gray] (4.5,0) to (4.5,1.3);
\draw [-, thin, gray] (5,0) to (5,1.3);
\draw [-, thin, gray] (5.5,0) to (5.5,1.3);
\draw [-, thin, gray] (6,0) to (6,1.3);
\draw [-, thin, gray] (6.5,0) to (6.5,1.3);
\draw [-, thin, gray] (7,0) to (7,1.3);
\draw [-, thin, gray] (7.5,0) to (7.5,1.3);
\draw [-, thin, gray] (8,0) to (8,1.3);
\draw [-, thin, gray] (8.5,0) to (8.5,1.3);
\draw [-, thin, gray] (9,0) to (9,1.3);
\draw [-, thin, gray] (9.5,0) to (9.5,1.3);
\draw [-, thin, gray] (10,0) to (10,1.3);
\draw [-, thick] (0,0) to (4,0) to (4,0.5) to (7, 0.5) to (7, 1) to (10,1);
\fill (0,0) circle(0.08);
\fill (10,1) circle(0.08);
\node at (10.8, 1.2)  {$(n,k)$};
\node at (0, -0.3)  {$0$};
\node at (4, -0.3)  {$t_1n$};
\node at (7, -0.3)  {$t_2n$};
\node at (10, -0.3)  {$n$};
\end{tikzpicture}
\caption{Thin DLPP}
\label{fig:thin}
\end{figure}

For the case of fixed $k$ and large $n$, a directed path looks like the one in Figure \ref{fig:thin}. Since the sum of $w_s$ on each row converges to a Brownian motion, Donsker's theorem implies that 
\beqq
	 \frac{L(n,k)}{\sqrt{n}} \Rightarrow   D_k \quad \text{where $D_k= \sup_{0=t_0\le t_1\le \cdots\le t_k=1} \sum_{i=1}^k (B_i(t_i)- B_i(t_{i-1}))$}
\eeqq
and $B_i(t)$, $i=1, \cdots, k$, are independent Brownian motions. 
On the other hand, it is known \cite{Baryshnikov01} that  
\beqq
	\lim_{k\to \infty} \Prob( (D_k- 2\sqrt{k})k^{1/6} \le x) = \FTW(x). 
\eeqq
This limit is a consequence of an explicit formula of the exponential DLPP by letting $n\to\infty$ first and then taking $k\to \infty$. 
We prove Theorem \ref{thm:thin} by showing that if $n, k\to \infty$ but $k$ grows slowly enough, we can take $n\to \infty$ first and $k\to \infty$ later. 
This argument is achieved by the Skorohod embedding or the Koml\'os--Major--Tusn\'ady embedding. 
However, the proof breaks down if $\alpha\ge 3/7$ since we cannot ignore the upward parts of the paths anymore. 
For $k=O(n)$, limit theorems are proved only for a few examples. 

\subsection{Interacting particle systems}

The particles in the TASEP move only one site to the right. Consider a more general finite-range exclusion process in which a particle at site $0$ can potentially move to site $v$ at rate $p(v)$. 
Assume that $\{v: p(v)>0\}$ is a finite set  generating $\Z$ additively, and $\sum_v v p(v)\neq 0$. 
In a recent paper \cite{Quastel-Sakar20}, Quastel and Sakar proved a KPZ limit theorem for finite-range exclusion processes started from a certain class of initial conditions. 
They compared the transition probabilities of the general process with those of the TASEP using energy estimates. 
This work is the first universality result in the $k=O(n)$ regime. It is exciting to see how the method further generalizes.  

\subsection{TASEP, Coulomb gas, and random matrices} \label{sec:TASEPRM}

In proving Theorem \ref{thm:Johansson}, Johansson also proved an unexpected connection of TASEP to Coulomb gas and random matrices \cite{Johansson00}. 
Consider the probability density function on $\R_+^n$ given by 
\beq \label{eq:Coulomb}
	p(x_1, \cdots, x_n) = c_{n, m} e^{2\sum_{1\le i<j\le n} \log |x_j-x_i|- \sum_{i=1}^n V(x_i) }, 
	\quad V(x)= x-(m-n)\log x, 
\eeq
where $c_{m,n}$ is the normalization constant.  
Let $x_{\max}=\max\{x_1, \cdots, x_n\}$. 
Johansson proved that for $m\ge n$, the last passage time $L(m,n)$ of the exponential DLPP has the same distribution as $x_{\max}$. 
The density function \eqref{eq:Coulomb} is said to define a  Coulomb gas with potential $V$ on $\R_+$ since the term $\log |x_i-x_j|$ is the 2d Coulomb potential of two equal charges at $x_i$ and $x_j$. 
The function $V(x)$ represents the confining potential. 

Let $X$ be an $n\times m$ random matrix with entries that are independent complex normal variables of mean zero and variance $1/2$. 
The random matrix $W=XX^*$ is called the complex Wishart matrix, and its eigenvalue density function is precisely \eqref{eq:Coulomb} \cite{Forrester10}. 
Thus, $L(m,n)$ has the same distribution as the largest eigenvalue of a complex Wishart matrix. 
See also \cite{Borodin-Peche08}. 
The connection of the exponential DLPP to Coulomb gas and the random matrix is special, and we do not expect to hold for general random variables $w_s$. 

There are universality results for both Coulomb gas and random matrices. 
The Tracy-Widom limit theorem is proved for the Coulomb gases with a general potential $V$. 
The paper \cite{Deift-Kriecherbauer-McLaughlin-Venakides-Zhou99} showed the limit theorem for generic analytic potentials, and \cite{Baik-Kriecherbauer-McLaughlin-Miller07} proved for discrete Coulomb gases in which particles are restricted to be only on a discrete set. 

Universality is a central question in random matrix theory, and there have been remarkable successes. 
The largest eigenvalue of a large class of random Hermitian matrices with independent entries converges to the Tracy-Widom distribution  \cite{Soshnikov99, Tao-Vu10, Erdos-Knowles-Yau-Yin12, Lee-Yin14}.

\subsection{Universality in many directions} \label{sec:univs}

We discussed that the TASEP with the step initial condition is connected to several areas: 
\begin{itemize}
\item interacting particle system
\item 1+1 random growth process
\item two-dimensional directed last passage percolation and directed polymer
\item Coulomb gas
\item random matrix
\end{itemize}
The TASEP also has interpretations as a random tiling model, and non-intersecting paths \cite{Johansson02, Gorin21}. 
We expect universality results to hold in all of these seven areas. 
The meaning of universality is  different in each area. For example, in random matrix theory, the largest eigenvalue of any random Hermitian matrix with independent and identically distributed entries with 4+$\epsilon$ finite moments converges to the same limit, the Tracy-Widom distribution. 
On the other hand, the 2d field limit of interacting particle systems depends on the initial condition, a special case of which has the Tracy-Widom distribution as its marginal. 

Even though many universality results are proved for Coulomb gases and random matrices, it remains to establish similar results for other areas and develop a general theory that encompasses all of these areas and possibly more.



\subsection*{Acknowledgements}

The author would like to thank Percy Deift, Zhipeng Liu, Andrei Prokhorov, and Guilherme Silva for helpful conversations. 
This work was partially supported by the NSF grant DMS-1954790.


\bibliographystyle{emss}


\def\cydot{\leavevmode\raise.4ex\hbox{.}}








\end{document}